\documentclass[12pt]{amsart}
\usepackage{amsfonts}
\usepackage{amssymb}

\usepackage{color}
\newtheorem{thm}{Theorem}[section]

\newtheorem{lem}[thm]{Lemma}
\newtheorem{prop}[thm]{Proposition}
\theoremstyle{definition}
\newtheorem{defn}[thm]{Definition}
\theoremstyle{remark}

\newcommand{\R}{\mathbb R}

\renewcommand{\span}{\mathrm{span}}
\newcommand{\X}{\mathfrak X}

\begin{document}

\title[A Classification of Riemannian QC-manifolds]
{A Classification of Riemannian manifolds of quasi-constant sectional curvatures}%
\author{Georgi Ganchev and Vesselka Mihova}%
\address{Bulgarian Academy of Sciences, Institute of Mathematics and
Informatics, Acad. G. Bonchev Str. bl. 8, 1113 Sofia, Bulgaria}%
\email{ganchev@math.bas.bg}%
\address{Faculty of Mathematics and Informatics, University of Sofia,.
J. Bouchier Str. 5, (1164) Sofia, Bulgaria}
\email{mihova@fmi.uni-sofia.bg}

% -------------------------------------------------------------------------

\begin{abstract}
Riemannian manifolds of quasi-constant sectional curvatures (QC-manifolds) are
divided into two basic classes: with positive or negative horizontal sectional
curvatures. We prove that the Riemannian QC-manifolds with positive horizontal
sectional curvatures are locally equivalent to canal hypersurfaces in Euclidean
space, while the Riemannian QC-manifolds with negative horizontal sectional
curvatures are locally equivalent to canal space-like hypersurfaces in Minkowski
space. We prove that the local theory of conformally flat Riemannian manifolds,
which can be locally isometrically embedded as hypersurfaces in Euclidean or
Minkowski space, is equivalent to the local theory of Riemannian QC-manifolds.
These results give a local geometric classification of conformally flat
hypersurfaces in Euclidean space and conformally flat space-like hypersurfaces
in Minkowski space.

\end{abstract}
\maketitle
% ----------------------------------------------------------------------------
\thispagestyle{empty}
\section{Introduction}
Conformally flat $n$-dimensional Riemannian manifolds appear as hypersurfaces
in two standard models of flat spaces: Euclidean or Minkowski space. Historically,
there were many attempts to describe conformally flat hypersurfaces, especially in
Euclidean space. Essential steps in this direction were made by Cartan \cite{C},
Schouten \cite{Sch1}. Kulkarni in \cite{K} reached to a partial description
of conformally flat hypersurfaces in Euclideasn space dividing them into:
hypersurfaces of constant curvature; hypersurfaces of revolution; tubes. Yano and
Chen proved in \cite{CY1} that canal hypersurfaces in Euclidean space are special
conformally flat hypersurfaces. Compact conformally flat hypersurfaces in Euclidean
space were studied in \cite{K} and \cite{dC}.

In this paper we study the close relation between the local theory of Riemannian
manifolds of quasi-constant sectional curvatures and the local theory of conformally
flat Riemannian hypersurfaces in the Euclidean space $\R^{n+1}$ or in the Minkowski
space $\R^{n+1}_1$. We give a local classification of Riemannian manifolds of
quasi-constant sectional curvatures proving that they can locally be embedded as canal
hypersurfaces in $\R^{n+1}$ or $\R^{n+1}_1$. Thus we obtain a geometric description
of conformally flat hypersurfaces in Euclidean space and conformally flat space-like
hypersurfaces in Minkowski space.

Riemannian QC-manifolds are Riemannian manifolds $(M,g,\xi)$ endowed with a
unit vector field $\xi$ besides the metric $g$, satisfying the curvature
condition: the sectional curvatures at any point of the manifold only depend on
the point and the angle between the section and the vector $\xi$ at that point.
All tangent sections at a given point, which are perpendicular to the vector
$\xi$ at that point, have one and the same sectional curvature. We call
these sectional curvatures horizontal sectional curvatures.

Everywhere, in this paper, we consider the case $\dim M=n \geq 4$.

The structural group of Riemannian manifolds $(M,g,\xi)$ is $O(n-1)\times 1$
and two Riemannian manifolds $(M,g,\xi)$ and $(M',g',\xi')$ are equivalent if
there exists a diffeomorphism $f: \, M \rightarrow M'$ preserving both structures:
the metric $g$ and the vector field $\xi$. We call such a diffeomorphism
a $\xi$-isometry.

In \cite {GM} we proved the following statements:
\vskip 1mm
{\it Any canal hypersurface $M$ in the Euclidean space $\R^{n+1}$ is a Riemannian
QC-manifold with positive horizontal sectional curvatures.}
\vskip 1mm
{\it Any Riemannian QC-manifold with positive horizontal sectional curvatures
is locally $\xi$-isometric to a canal hypersurface in the Euclidean space $\R^{n+1}$.}
\vskip 1mm
{\it The first problem} we treat here is to give a local classification of
Riemannian $QC$-manifolds with negative horizontal sectional curvatures.

In Section 3 we introduce canal space-like hypersurfaces in the Minkowski
space $\R_1^{n+1}$ and divide them into three types. In Subsections 3.1 - 3.3
we study these three types of canal space-like hypersurfaces and show that
\vskip 1mm
{\it Any canal space-like hypersurface in the Minkowski space $\R^{n+1}_1$ is a
Riemannian QC-manifold with negative horizontal sectional curvatures.}
\vskip 1mm
The basic results are proved in Section 4. The local classification  of Riemannian
$QC$-manifolds with negative horizontal sectional curvatures is given by Theorem 4.1:
\vskip 1mm
{\it Any Riemannian QC-manifold with negative horizontal sectional curvatures
is locally $\xi$-isometric to a canal space-like hypersurface in the Minkowski
space $\R^{n+1}_1$.}
\vskip 1mm
{\it The second problem} we deal with is to obtain a geometric
description of conformally flat hypersurfaces in Euclidean space and conformally
flat space-like hypersurfaces in Minkowski space. Using results of Cartan and
Schouten, we obtain the second fundamental form of the hypersurface into
consideration. This allows us to give a local geometric classification of
conformally flat Riemannian hypersurfaces in Euclidean or Minkowski space:
\vskip 1mm
{\it Any conformally flat hypersurface in Euclidean space, which is free of
umbilical points, locally is a part of a canal hypersurface.

Any conformally flat space-like hypersurface in Minkowski space, which is free
of  umbilical points, locally is a part of a canal space-like hypersurface.}

The picture of the local isometric embeddings of a conformally flat Riemannian
manifold into Euclidean or Minkowski space can be described briefly as follows.

Let $(M,g)$ be a conformally flat Riemannian manifold, free of points in which
all sectional curvatures are constant. The manifold $(M,g)$ can be locally
isometrically embedded into $\R^{n+1}$ ($\R^{n+1}_1$) if and only if its Ricci
operator has two different from zero eigenvalues at every point: one of them
of multiplicity $n-1$, the other of multiplicity $1$. The latter eigenvalue
generates a unit vector field $\xi$, such that $(M,g,\xi)$ is a Riemannian
QC-manifold with positive (negative) horizontal sectional curvatures. Any two
isometrical realizations of $(M,g)$ are locally congruent.

Generalizing, we obtain that a conformally flat Riemannian manifold is locally
isometric to a hypersurface into $\R^{n+1}$ ($\R^{n+1}_1$) if and only if its
Ricci operator at any point has a root of multiplicity at least $n-1$. This fact
gives an approach to further investigations of conformally flat Riemannian
manifolds studying the spectrum of their Ricci operator.

It is interesting to mention Riemannian subprojective manifolds forming a subclass
of Riemannian QC-manifolds characterized by the condition: the structural vector
field $\xi$ is geodesic. If $(M,g,\xi)$ is a Riemannian subprojective manifold
with scalar curvature $\tau \, (d\tau \neq 0)$, then the structural vector field
$\xi$ is collinear with ${\rm grad} \, \tau$. Any Riemannian subprojective manifold
is locally isometric (up to a motion) to a rotational hypersurface in Euclidean
space or in Minkowski space.

\section{Preliminaries}

Let $(M,g,\xi) \, (\dim M = n \geq 4)$ be a Riemannian manifold with metric $g$ and
a unit vec-tor field $\xi$. The structural group of these manifolds is
$O(n-1)\times 1$. $T_pM$ and $\X M$ will stand for the tangent space to $M$
at a point $p$ and the algebra of smooth vector fields on $M$, respectively.
The 1-form corresponding to the unit vector $\xi$ is denoted by $\eta$, i.e.
$\eta(X)=g(\xi,X), \; X \in \X M$. The distribution of the 1-form $\eta$ is
denoted by $\Delta$, i.e.
$$\Delta(p)=\{X \in T_pM:\,\eta(X)=0\}.$$

The orthogonal projection of a vector field $X \in \X M$ onto the distribution
$\Delta$ is denoted by the corresponding small letter $x$, i.e.
$$X=x+\eta(X)\,\xi.\leqno(2.1)$$

Any section $E$ in $T_pM$ determines an angle $\angle (E,\xi)$. Then the notion
analogous to the notion of a Riemannian manifold of constant sectional
curvatures is described as follows \cite{GM}.

\begin{defn} A Riemannian manifold $(M,g,\xi)$ $(\dim M \geq 3)$ is said to be of
quasi-constant sectional curvatures (a Riemannian $QC$-manifold) if for an arbitrary
2-plane $E$ in $T_pM,\, p\in M$, with $\angle (E,\xi)= \varphi$, the sectional curvature
of $E$ only depends on the point $p$ and the angle $\varphi$.
\end{defn}

Let $\nabla$ be the Levi-Civita connection of the metric $g$ and $\mathcal R$
be its Riemannian curvature tensor. The structure $(g, \xi)$ generates the
following tensors $\pi$ and $\Phi$:
$$\begin{array}{ll}
\pi(X,Y,Z,U)=&g(Y,Z)g(X,U)-g(X,Z)g(Y,U),\\
[3mm]
\Phi(X,Y,Z,U)=&g(Y,Z)\eta(X)\eta(U)-g(X,Z)\eta(Y)\eta(U)\\
[3mm]
 & + \,g(X,U)\eta(Y)\eta(Z)-g(Y,U)\eta(X)\eta(Z);
\quad X,Y,Z,U \in \X M.\end{array}$$
These tensors have the symmetries of the curvature tensor $\mathcal R$ and are invariant
under the action of the structural group of the manifold.

Riemannian manifolds of quasi-constant sectional curvatures are characterized by the
following statement \cite{GM}:

\begin{prop}
A Riemannian manifold $(M,g,\xi)$ is of quasi-constant sectional curvatures if and
only if its curvature tensor has the form
$$\mathcal R=a\,\pi + b \, \Phi,\leqno(2.2)$$
where $a$ and $b$ are some functions on $M$.
\end{prop}

Let $(M,g,\xi)$ $(\dim M=n\geq 4)$ be a Riemannian manifold of quasi-constant
sectional curvatures. This means that the curvature tensor $\mathcal R$ of $g$
has the form (2.2). If $b \neq 0$, then the manifold $(M,g,\xi)$ has the
properties \cite{GM}:
\begin{itemize}
\item The distribution of the function $a$ is the structural distribution $\Delta$:
$$da=\xi(a)\,\eta.\leqno(2.3)$$
\par
\item The distribution $\Delta$ is involutive, i.e.
$$d\eta(x,y)=0, \quad x,y \in \Delta.\leqno(2.4)$$
\par
\item If $\theta$ is the 1-form defined by $\theta(X)=d\eta(\xi, X), \; X \in \mathcal XM$,
then $d\eta= \theta \wedge \eta$ and
$$\theta(x) = d\eta(\xi, x)=\frac{1}{b}\,d b(x),\quad x\in \Delta.\leqno (2.5)$$
\par
\item The integral submanifolds of the distribution $\Delta$ are totally umbilical
in $M$, i.e.
$$\nabla_x\xi=k \,x, \quad k=\frac{\xi(a)}{2b}, \quad x \in \Delta.\leqno(2.6)$$
\par
\item The distribution of the function $k$ is the structural distribution $\Delta$:
$$dk=\xi(k)\,\eta.\leqno(2.7)$$
\end{itemize}

Let $S_p$ be the maximal integral submanifold of the distribution $\Delta$,
containing a given point $p \in M$, and $\mathcal K$ be the curvature tensor of
the Riemannian manifold $(S_p,g)$. Then we have:

(i) All sections tangent to $S_p$ have one and the same sectional curvature $a(p)$
with respect to the tensor $\mathcal R$. We say that the function $a(p)$ is
the horizontal sectional curvature of the manifold.

(ii) All sections tangent to $S_p$ have one and the same sectional
curvatures $a(p)+k^2(p)$ with respect to the tensor $\mathcal K$.

Proposition 2.2 implies the following statement.

\begin{prop}
A Riemannian QC-manifold $(M,g,\xi)$, free of points in which the sectional
curvatures are constant, i.e. $b \neq 0$, is characterized by the following
two conditions:

- $(M,g)$ is conformally flat;
\vskip 2mm
- the Ricci operator $\rho$ of $(M,g)$ at any point has two non-zero roots:

$(n-1)a+b$, \, of multiplicity \, $n-1$, \, which generates the distribution $\Delta$:
$$\rho(x)= [(n-1)a+b]\,x, \quad x \in \Delta;$$

$(n-1)(a+b)$, \, of multiplicity \, $1$, which generates the structural vector field $\xi$:
$$\rho(\xi)=(n-1)(a+b)\,\xi.$$
\end{prop}

Proposition 2.3 implies that the notion of a QC-manifold is a notion in Riemannian
geometry. The next statement is an immediate consequence from this proposition.

\begin{thm} Let $(M,g,\xi)$ and $(\bar M, \bar g, \bar \xi)$ be two QC-manifolds
free of points in which the sectional curvatures are constant. If
$\varphi : \, M \; \rightarrow \; \bar M$ is an isometry, then it is a $\xi$-isometry,
i.e. $\varphi_* \xi=\bar \xi$.
\end{thm}

The above mentioned geometric functions $a$ and $a+k^2$ on $(M,g,\xi)$ generate
four basic classes of Riemannian manifolds of quasi-constant sectional curvatures
characterized by the conditions:
\begin{itemize}
\item[1)] \; $a > 0$;

\item[2)] \; $a< 0, \quad a+k^2>0$;

\item[3)] \; $ a+k^2<0$;

\item[4)] \; $a+k^2=0$.
\end{itemize}
\vskip 1mm

The class of Riemannian QC-manifolds contains the remarkable subclass of Riemannian
subprojective manifolds. V. Kagan \cite{K1, K2} called an $n$-dimensional  space
$A_n$ with symmetric linear connection $\nabla$ a subprojective space if there
exists locally a coordinate system with respect to which every geodesic of
$\nabla$ can be represented by $n-2$ linear equations  and another equation,
that need not be linear (see also \cite{Sch2}). P. Rachevsky \cite{R} proved
necessary and sufficient conditions characterizing Riemannian subprojective spaces.
T. Adati \cite{A} studied Riemannian subprojective manifolds concerning
concircular and torse-forming vector fields.

As Riemannian QC-manifolds $(M,g,\xi)$ the Riemannian subprojective manifolds are
characterized by any of the following additional properties \cite{GM}:
\begin{itemize}
\item[i)] \; $db = \xi(b)\,\eta$;

\item[ii)] \; the vector field $\xi$ is geodesic (on $M$);

\item[iii)] \; the 1-form $\eta$ is closed.
\end{itemize}

Let $\tau$ be the scalar curvature of a Riemannian subprojective manifold.
If $d \tau \neq 0$, then the structural distribution $\Delta$ is the distribution
of the 1-form $d \tau$ and the vector field ${\rm grad}\,\tau$ is an eigenvector
of the Ricci operator at every point.

\section{Canal space-like  hypersurfaces in Minkowski space}

A hypersurface $M \,(\dim M=n)$ in the Minkowski space $\mathbb{R}^{n+1}_1$ is said to be
{\it space-like} (or \emph{Riemannian}) if the induced metric on $M$ is positive definite.
The normal vector field to a space-like hypersurface $M$ in the Minkowski space
$\mathbb{R}^{n+1}_1$ is necessarily time-like.

In this section we study the envelope of a one-parameter family of space-like hyperspheres
$\{S^n(s)\},\; s \in J\subset \mathbb{R}$ in $\R^{n+1}_1$, given as follows
$$S^n(s):\quad (Z - z(s))^2 = -R^{2}(s), \quad R(s)>0,$$
where $z=z(s)$ is the center and $R(s)$ is the radius of the corresponding hypersphere
$S^n(s)$.

Let the cross-section of a space-like hypersphere $S^n$ with a hyperplane in the Minkowski
space $\mathbb{R}^{n+1}_1$ be an $(n-1)$-dimensional surface. We have:
\begin{itemize}
\item[1)]\, The cross-section of a space-like hypersphere $S^n$ with a space-like
hyperplane $R^n$ is a Euclidean hypersphere $S^{n-1}$ in $\R^n$, and $S^{n-1}$ is
of positive constant sectional curvatures.

\item[2)]\, The cross-section of a space-like hypersphere $S^n$ with a time-like
hyperplane $R^n_1$ is a hyperbolic hypersphere $H^{n-1}$ in $\R^n_1$, and
$H^{n-1}$ is of negative constant sectional curvatures.

\item[3)]\, The cross-section of a space-like hypersphere $S^n$ with a light-like
hyperplane $R^n_0$ is a parabolic hypersphere $P^{n-1}$ in $\R^n_0$, and $P^{n-1}$
is of zero sectional curvatures.
\end{itemize}

We shall describe in more details the cross-section $P^{n-1}$ of a space-like
hypersphere $S^n(z,R)$ with a light-like hyperplane $\R^n_0$. It is clear that
$\R^n_0$ can not pass through the center $z$ of $S^n$. The pair $(\R^n_0, g)$
is an $n$-dimensional affine space with metric $g$, whose rank equals $n-1$.
This means that $\R^n_0$ contains a light-like direction $U$, determined by a
given light-like vector $t$. The light-like direction $U$ can also be considered
as a point at infinity in the infinite hyperplane of $\R^n_0$. Any hyperplane
$E^{n-1}$ of $\R^n_0$, which does not contain $U$, is a Euclidean hyperplane,
i.e. it can be endowed with a basis $e_1, ..., e_{n-1}$, satisfying the property
$g(e_i, e_j)=\delta_{ij}, \, i,j=1, ..., n-1$, \, $\delta_{ij}$ being the
Kronecker's deltas.

Let $E^{n-1}$ be a Euclidean hyperplane in $\R^n_0$ with a fixed point $T \in E^{n-1}$
and an orthonormal basis $e_1,...,e_{n-1}$. Adding the light-like vector $t$, we obtain
a coordinate system $Te_1,...,e_{n-1}t$ in $\R^n_0$. If $Z(z_1,...,z_{n-1};z_n)$ is the
position vector of any point $Z$ in $\R^n_0$, then we consider
the quadrics $P^{n-1}(q)$ in $\R^n_0$, given by the equation
$$ P^{n-1}(q):\; z_1^2+...+z_{n-1}^2-2q\,z_n=0, \quad q={\rm const} > 0.$$

The one-parameter family of quadrics $P^{n-1}(q)$ is characterized by the properties:
\begin{itemize}
\item[(i)] $P^{n-1}(q)$ is a quadric, which is tangent to the infinite hyperplane of $\R^n_0$
at $U$ and to the hyperplane $E^{n-1}$ at $T$;

\item[(ii)] The cross-section of $P^{n-1}(q)$ with any Euclidean hyperplane
$z_n= {\rm const} > 0$ (parallel to $E^{n-1}$) is a Euclidean hypersphere
in this hyperplane.
\end{itemize}

We call these quadrics \emph{parabolic hyperspheres} of the light-like hyperplane
$(\R^n_0, g)$.

The parabolic hyperspheres  have the following remarkable property:

\begin{prop} Any parabolic hypersphere in a light-like hyperplane $\R^n_0$ is a flat
$(n-1)$-dimensional Riemannian manifold.
\end{prop}

\emph{Proof:} Since the only tangent hyperplane to the parabolic hypersphere
$P^{n-1}(q)$, which contains $U$, is the infinite hyperplane of $\R^n_0$, then
$(P^{n-1}(q), g)$ is an $(n-1)$-dimensional Riemannian manifold.

We consider the projection
$$\pi: \quad P^{n-1}(q) \; \rightarrow \, E^{n-1}$$
of the parabolic hypersphere onto the Euclidean hyperplane $E^{n-1}$, parallel
to the direction $U$. It is an easy verification that the projection $\pi$ is
an isometry between the Riemannian manifolds $(P^{n-1}(q), g)$ and $(E^{n-1}, g)$,
excluding the common point $T$. This implies the assertion. \qed
\vskip 2mm
Next, we call the $(n-1)$-dimensional cross-sections of a space-like hypersphere with a
hyperplane \emph{spheres of codimension two} and use the common denotation $S^{n-1}$.

Let $M=\{S^{n-1}(s)\},\; s \in J\subset \mathbb{R}$ be a space-like hypersurface in
$\mathbb{R}^{n+1}_1$, which is a one-parameter family of spheres $S^{n-1}(s)$ of
codimension two. Any sphere $S^{n-1}(s)$ is said to be {\it a spherical generator} of $M$.

At first canal surfaces in $\R^3$ have been introduced and studied in the classical
works of Enneper \cite{E1, E2, E3, E4}. We use the following definition:
\begin{defn}
A space-like hypersurface $M=\{S^{n-1}(s)\},\; s \in J\subset \mathbb{R}$ in
$\mathbb{R}^{n+1}_1$ is said to be {\it a canal space-like hypersurface} if the normals
to $M$ at the points of any fixed spherical generator pass through a fixed point.
\end{defn}

Let now $ Z=Z(s; u^1, u^2,...,u^{n-1}),\;s \in J,\; (u^1, u^2,...,u^{n-1})\in D$ be
the position vector field of a canal space-like hypersurface $M$. The partial
derivatives of $Z$ are denoted as follows:
$Z_s=\frac{\partial Z}{\partial s},\; Z_i=\frac{\partial Z}{\partial u^i}; \; i=1,...,n-1,$
and similar denotations are used for other vector functions.

Denoting by $z(s), \, s\in J$ the common point of the normals to $M$ at the points of
any spherical generator $S^{n-1}(s)$,
we consider the space-like hypersphere $S^n(s)$ with center $z(s)$ containing $S^{n-1}(s)$.
If $R(s)$ is the radius of $S^n(s)$, then the position vector $Z$ of $M$ satisfies
the equality
$$ (Z - z(s))^2 = -R^{2}(s), \quad R(s)>0, \quad s \in J\subset \mathbb{R}.\leqno (3.1)$$
Differentiating (3.1) with respect to the parameter $s$, we get
$$(Z-z(s))Z_s - (Z-z(s))z'(s)=-R(s)R'(s). \leqno (3.2)$$
Under the condition that the normal to $M$ at any point of a fixed generator $S^{n-1}(s)$
is collinear with $Z-z(s)$, the equalities (3.1) and (3.2) are equivalent to
$$\begin{array}{l}
(Z-z(s))^2=-R^2(s),\\
[3mm]
(Z-z(s))z'(s)=R(s)R'(s).
\end{array} \leqno(3.3)$$

A space-like hypersurface $M$ in $\mathbb{R}^{n+1}_1$ is said to be the
{\it envelope} of a one-parameter family of space-like hyperspheres
$\{S^n(z(s),R(s))\}, \, s \in J$ if the position vector
$Z(s;u^1,...,u^{n-1})$ of $M$ satisfies the equations (3.3).

Let $M$ be a space-like hypersurface, which is the envelope of a one-parameter
family of space-like hyperspheres $\{S^n(z(s),R(s))\}, \, s \in J$. It follows
from (3.3) that $M$ is a one parameter family of spheres $S^{n-1}(s), \, s\in J$.
Differentiating the first equality of (3.3), we have
$$(Z-z) Z_s = 0,\quad (Z-z) Z_i=0, \; i=1,..,{n-1},$$
which shows that the time-like vector field $Z-z$ at the points of any generator
$S^{n-1}(s)$ of $M$ is normal to both: the hypersurface $M$ and
the hypersphere $S^n(s)$.

Hence,  as in the classical case \cite{E1, E2, E3, E4, vL}, we have
\begin{lem}
A space-like hypersurface $M$ in $\mathbb{R}^{n+1}_1$ is canal if and only if
it is the envelope of a one-parameter family of space-like hyperspheres.
\end{lem}

Let $M$ be a space-like canal hypersurface, given by (3.3). We denote the tangent
vector to the curve of centers $z(s)$ as usual by $t(s)=z'(s)$. The unit  normal
vector field $N$  to $M$ is collinear with $Z-z$ and we always choose
$$N = - \frac{Z-z}{R}. \leqno {(3.4)}$$

In view of (3.3), the vector field $N$ has the properties:
$$N^2=-1,\quad Nt=-R'.$$

Differentiating (3.4), we have
$$\begin{array}{l}
\displaystyle{N_i=-\frac{1}{R}\, Z_i, \; i=1,...,n-1, }\\
[4mm]
Z_s+R\, N_s=t-R' N.
\end{array}\leqno (3.5)$$

The second equality in (3.5) means that the vector field \,$t-R' N$
is tangent to $M$. Since the normals to $M$ at the points of a
spherical generator cannot be parallel to the vector $t$, then the
vector field \,$t-R' N$ is space-like and $(t-R'\,N)^2> 0$.
Furthermore, the second equality in (3.3) implies that $t\,Z_i=0, \;
i=1,...,n-1$, and therefore \,$t-R' N$ is perpendicular to all
$Z_i$.

We introduce the unit tangent vector field $\xi$  as follows:
$$ \xi := \frac{1}{\sqrt{(t-R'\,N)^2}}\,(t-R'\,N).  \leqno{(3.6)}$$

Then the distribution $\Delta:=\{x\in T_p M \,: \,x \perp \xi\}$
is exactly \; $\Delta = \span\{Z_1,...,Z_{n-1}\}$.

For the purposes of our investigations we need to introduce three types of canal
space-like hypersurfaces.

\begin{defn}
A  canal space-like hypersurface $M$ in $\mathbb{R}^{n+1}_1$, given by (3.3),
is said to be {\it a canal space-like hypersurface of elliptic, hyperbolic} or
{\it parabolic type} if the curve $z=z(s), \, s\in J$ of the centers of the
hyperspheres is time-like, space-like or light-like, respectively.
\end{defn}

Rotational space-like hypersurfaces are introduced in a natural way:

\begin{defn}
A  canal space-like hypersurface $M$ in $\mathbb{R}^{n+1}_1$, given by (3.3),
is said to be {\it a rotational space-like hypersurface} if the curve
$z=z(s), \, s\in J$ of the centers of the hyperspheres lies on a straight line.
\end{defn}

Any of the three types of canal space-like hypersurfaces generates the corresponding
subclass of rotational space-like hypersurfaces.

\subsection{Canal space-like hypersurfaces of elliptic type in Minkowski space}
Let $M$ be a canal space-like hypersurface in $\mathbb{R}^{n+1}_1$ of elliptic type,
given by (3.3). The curve of centers $z=z(s),\;s \in J$, parameterized by its
natural parameter, satisfies the condition $z'^2= t^2 = -1$.

Since $(t-R'\,N)^2=R'^2-1 > 0$, then the function $R(s)$ in the case of a canal
space-like hypersurface of elliptic type satisfies the inequalities
$$R^2(s) > 0, \quad R'^2(s)-1 > 0; \quad s \in J.$$

Next we find the second fundamental form of $M$.

Let $\nabla '$ be the standard flat Levi-Civita connection in ${\bf \mathbb{R}^{n+1}_1}$
and $h$ be the second fundamental tensor of $M$. The Levi-Civita connection of the
induced metric on the hypersurface $M$ is denoted by $\nabla$. Taking into account
(3.5) and (3.6), we get
$$\nabla'_{Z_i} N=N_i=-\frac{1}{R} \, Z_i,\quad \nabla'_{Z_i} \xi =\xi_i=
-\frac{R'}{\sqrt{R'^2-1}} \, N_i= \frac{R'}{R\sqrt{R'^2-1}} \, Z_i.$$
These equalities can be written as follows:
$$\nabla '_x N=- \frac{1}{R} \,x,\quad \nabla'_x \xi =\nabla_x \xi
=\frac{R'}{R\sqrt{R'^2-1}}\,x\, =k\,x,\quad x \in \Delta,$$
and the function $k$ is
$$k=\frac{R'}{R\sqrt{R'^2-1}} \,. \leqno(3.1.1)$$
Hence, the shape operator $A$ of $M$ satisfies
$$Ax=\frac{1}{R} \, x, \quad x \in \Delta. \leqno (3.1.2)$$

Since $g(N,N)=-1$, then
$$h(x,y) = -g(Ax,y)=-\frac{1}{R}\,g(x,y), \quad x,y \in \Delta. \leqno(3.1.3)$$

The equality (3.1.2) means that the tangent space $\Delta$ is invariant with
respect to the shape operator $A$. This implies that the vector field $\xi$
is also an eigenvector field of $A$, i.e.
$$A\xi= \nu\,\xi. \leqno (3.1.4)$$

Assuming the standard summation convention, we can put
$$\xi=\phi^i Z_i +\phi Z_s, \quad \phi \neq 0 \leqno(3.1.5)$$
for some functions $\phi^1,..., \phi^{n-1}; \, \phi$\,
on $M$. Since $\xi$ is perpendicular to all $Z_i$, we have
$$\phi \, (\xi Z_s)=1. \leqno(3.1.6)$$

Taking into account (3.1.5), we compute
$$\nabla'_{\xi} N= \phi^i N_i +\phi N_s
=-\frac{1}{R} \, \phi^i Z_i + \phi N_s.
\leqno (3.1.7)$$

On the other hand, because of (3.5) and (3.6), we have
$$Z_s + R\,N_s= \sqrt{R'^2-1}\,\xi.
\leqno (3.1.8)$$
In view of (3.1.5) and (3.1.8) equality (3.1.7) implies that
$$\nabla'_{\xi} N= -\frac{1}{R}(1-\phi \sqrt{R'^2-1})\xi=
-\nu\, \xi,$$
and
$$\nu - \frac{1}{R}= -\frac{\sqrt{R'^2-1}}{R} \,\phi \,. \leqno(3.1.9)$$

Using (3.1.2), (3.1.3) and (3.1.4), we obtain the shape operator of $M$:
$$AX=\frac{1}{R}\,X+\left(\nu-\frac{1}{R}\right)\eta(X)\,\xi,\quad X\in \mathcal{X}M.$$

The last equality and (3.1.9) imply that the second fundamental tensor of $M$
has the form
$$h(X,Y)=-\frac{1}{R}\,g(X,Y)+\phi\,\frac{\sqrt{R'^2-1}}{R} \, \eta(X)\eta(Y),\quad
X, Y\in \mathcal{X}M. \leqno (3.1.10)$$

Further we replace (3.1.10) into the Gauss equation for the hypersurface $M$, and
taking into account (3.1.6), we obtain the curvature tensor $\mathcal R$ of
a canal space-like hypersurface $M$ of elliptic type:
$$\mathcal{R}=-\frac{1}{R^2}\,\pi+\frac{\sqrt{R'^2-1}}{R^2(\xi Z_s)}\,\Phi=a\pi+b\Phi.
\leqno (3.1.11)$$

Now (3.1.11) and (3.1.1) imply that
$$a=-\frac{1}{R^2} < 0, \qquad a+k^2=\frac{1}{R^2(R'^2-1)}>0.$$

Thus we obtained the following
\begin{prop}
Any canal space-like  hypersurface of elliptic type in $\mathbb{R}^{n+1}_1$ is
a Riemannian manifold of quasi-constant sectional curvatures with functions
$a<0$ and $a+k^2 >0$.
\end{prop}

Next we prove that the rotational space-like hypersurfaces of elliptic type are
Riemannian subprojective manifolds satisfying the conditions in
Proposition 3.6.

Using (3.1.8), we have
$$\xi Z_s+R \,(\xi N_s) = \sqrt{R'^2-1}.$$

In order to compute the function $\xi N_s$, we use the equality $\xi N_s + \xi_s N=0$.
Differentiating (3.6) by $s$, we find
$$\xi_s N=\frac{t'N +R''}{\sqrt{R'^2-1}}\,.$$
Therefore
$$\xi Z_s=\frac{RR''+R'^2-1+R(t'N)}{\sqrt{R'^2-1}}\,$$
and
$$b=\frac{R'^2-1}{R^2\{RR''+R'^2-1+R(t'N)\}}\,.$$

According to Proposition 3.6,
the hypersurface $M$ is a Riemannian QC-manifold.
Any Riemannian QC-manifold is subprojective if and only if the functions $a$ and $b$
generate one and the same distribution. Therefore, $M$ is subprojective if and only if
the function $b$ does not depend on the parameters $u^i; \, i=1,...,n-1$, i.e. $t'=0$.
Since $t'=0$ characterizes a straight line $c$, we obtain the following statement.
\begin{prop}
A canal space-like hypersurface $M$ of elliptic type in $\R_1^{n+1}$ is a
Riemannian subprojective manifold if and only if $M$ is a rotational space-like
hypersurface of elliptic type.
\end{prop}
Combining with Proposition 3.6, we have
\begin{prop}
Any rotational space-like hypersurface of elliptic type in $\mathbb{R}^{n+1}_1$ is
a subprojective Riemannian manifold  with functions $a<0$ and $a+k^2 >0$.
\end{prop}
The curvature tensor of a rotational space-like hypersurface of elliptic type has
the form
$$\mathcal R= -\frac{1}{R^2}\, \pi + \frac{R'^2-1}{R^2(RR''+R'^2-1)}\,\Phi.$$

\subsection{Canal space-like hypersurfaces of hyperbolic type in Minkowski space}
Let $M$ be a canal space-like hypersurface of hyperbolic type, given by (3.3).
The curve of centers $z=z(s),\;s \in J$, parameterized by its natural parameter,
satisfies the condition $z'^2= t^2 = 1$.

In the case considered, the inequality $(t-R'\,N)^2= R'^2+1 >0$ is always satisfied.
Hence, $R(s)$ satisfies the only condition $R(s)>0$.

We compute
$$\nabla '_x N=- \frac{1}{R}\,x,\quad \nabla'_x \xi =
\nabla_x \xi=\frac{R'}{R\sqrt{R'^2+1}}\, x = k\,x,
\quad x \in \Delta,$$
where the function $k(s)$  is
$$k=\frac{R'}{R\sqrt{R'^2+1}}.\,\leqno(3.2.1)$$
Therefore,
$$Ax=\frac{1}{R}\,x,\quad x \in \Delta\leqno (3.2.2)$$
and the vector field $\xi$ is an eigenvector for $A$:
$$A\xi= \nu\,\xi, \leqno (3.2.3)$$

Putting $\xi=\phi^i Z_i +\phi Z_s$, we compute
$$\nabla'_{\xi} N= -\frac{1}{R} \phi^i Z_i + \phi N_s, $$
and taking into account that
$$Z_s + R\,N_s= \sqrt{R'^2+1}\,\xi, $$
we find
$$\nabla'_{\xi} N= -\frac{1}{R}(1-\phi \sqrt{R'^2+1})\xi= -\nu\, \xi,$$
and
$$\nu=\frac{1}{R}(1-\phi \sqrt{R'^2+1}).$$

Using (3.2.2) and (3.2.3), we obtain the second fundamental form $h$ of the
hypersurface $M$:
$$h=-\frac{1}{R}\,g+\phi\frac{\sqrt{R'^2+1}}{R}\,\eta \otimes\eta.$$

Applying the Gauss equation and the equality $\phi \,(\xi Z_s)=1$, we calculate the
curvature tensor of the hypersurface $M$.
$$\mathcal{R}=-\frac{1}{R^2}\,\pi+\frac{\sqrt{R'^2+1}}{R^2\,(\xi Z_s)}\,\Phi=a\pi+b\Phi.
$$

Therefore, the function $a=-1/R^2$. In view of (3.2.1), we find
$$a+k^2=\frac{-1}{R^2(R'^2+1)}<0.$$

Thus we obtained the following statement.
\begin{prop}
Any canal space-like hypersurface of hyperbolic type in $\mathbb{R}^{n+1}_1$ is
a Riemannian manifold of quasi-constant sectional curvatures with function  $a+k^2 <0$.
\end{prop}

Next we prove that the rotational space-like hypersurfaces of hyperbolic type are
Riemannian subprojective manifolds satisfying the condition in Proposition 3.9.

Differentiating (3.6) with respect to $s$, we compute
$$\xi_s N=\frac{t'N+R''}{\sqrt{R'^2+1}}\,.\leqno(3.2.4)$$
Using the equality $\xi_s N + \xi N_s =0$, (3.2.4) and (3.5), we find
$$\xi Z_s=\frac{RR''+R'^2+1+R\,(t'N)}{\sqrt{R'^2+1}}$$
and
$$b=\frac{R'^2+1}{R^2\{RR''+R'^2+1+R\,(t'N)\}}.$$

Applying similar arguments as in Subsection 3.1, we conclude that $M$ is
subprojective if and only if $t'=0$, which characterizes a straight line $c$.

Thus we obtained the following statement.
\begin{prop}
A canal space-like hypersurface $M$ of hyperbolic type in $\R_1^{n+1}$ is a
Riemannian subprojective manifold if and only if $M$ is a rotational space-like
hypersurface of hyperbolic type.
\end{prop}
Combining with Proposition 3.9, we have
\begin{prop}
Any rotational space-like hypersurface of hyperbolic type in $\mathbb{R}^{n+1}_1$ is
a subprojective Riemannian manifold  with function $a+k^2 < 0$.
\end{prop}
The curvature tensor of a rotational space-like hypersurface of hyperbolic type has
the form
$$\mathcal R= -\frac{1}{R^2}\, \pi + \frac{R'^2+1}{R^2(RR''+R'^2+1)}\,\Phi.$$

\subsection{Canal space-like hypersurfaces of parabolic type in Minkowski space}
Let $M$ be a canal space-like hypersurface of hyperbolic type, given by (3.3).
The curve of centers $z=z(s),\;s \in J$, satisfies the condition $z'^2= t^2 =0$.

In this case  $(t-R'\,N)^2=R'^2>0$ and the function $R(s)$ satisfies the conditions
$R(s)>0$ and $R'(s)\neq 0$.

Next we find the second fundamental form of $M$.

We compute
$$\nabla '_x N=- \frac{1}{R}\,x,\quad \nabla'_x \xi = \nabla_x \xi=\frac{1}{R}\, x =k x,
\quad x \in \Delta $$
and find
$$Ax=\frac{1}{R}\,x,\quad x \in \Delta, \leqno (3.3.1)$$
$$k=\frac{1}{R} \leqno(3.3.2)$$
and
$$A\xi=\nu \, \xi. \leqno(3.3.3)$$

Further we again put
$$\xi=\phi^i Z_i +\phi Z_s$$
and compute
$$\nabla'_{\xi} N= \phi^i N_i +\phi N_s
=-\frac{1}{R} \phi^i Z_i + \phi N_s,\quad i=1, 2, ... ,n-1. \leqno (3.3.4)$$
Using the equality
$$Z_s= R'\,\xi - R\,N_s,$$
we obtain from (3.3.4) that
$$\nabla'_{\xi} N= -\frac{1}{R}(1-\phi R')\xi= -\nu\, \xi,$$
and
$$\nu=\frac{1}{R}(1-\phi R').\leqno (3.3.5)$$

Now equalities (3.3.1), (3.3.3) and (3.3.5) imply  that
$$h=-\frac{1}{R}\,g+\phi \frac {R'}{R} \,\eta\otimes\eta.$$

Finally, replacing $h$ into the Gauss equation and using the equality $\phi \, (\xi Z_s)=1$,
we find the curvature tensor of the hypersurface $M$ in the form
$$\mathcal{R}=-\frac{1}{R^2}\,\pi+\frac{R'}{R^2\,(\xi Z_s)} \,\Phi=a\pi+b\Phi,$$
which shows that $M$ is a Riemannian QC-manifold with function $a=-1/R^2$.
In view of (3.3.2) we find
$$a+k^2=0.$$

Thus we obtained the following statement.
\begin{prop}
Any canal space-like hypersurface of parabolic type in $\mathbb{R}^{n+1}_1$ is
a Riemannian manifold of quasi-constant sectional curvatures with function  $a+k^2 =0$.
\end{prop}

Next we prove that the rotational space-like hypersurfaces of hyperbolic type are
Riemannian subprojective manifolds satisfying the condition in Proposition 3.12.

Differentiating (3.6) with respect to $s$, we get
$$\xi_s N=\frac{t'N+R''}{R'}\,.\leqno(3.3.6)$$
Using the equality $\xi_s N + \xi N_s =0$, (3.3.6) and (3.5), we find
$$\xi Z_s=\frac{RR''+R'^2+R'\,(t'N)}{R'}$$
and
$$b=\frac{R'^2}{R^2\{RR''+R'^2+R'\,(t'N)\}}.$$

Applying similar arguments as in Subsection 3.1, we conclude that $M$ is
subprojective if and only if $t'=0$, which characterizes a straight line $c$.

Thus we obtained the following statement.
\begin{prop}
A canal space-like hypersurface $M$ of parabolic type in $\R_1^{n+1}$ is a
Riemannian subprojective manifold if and only if $M$ is a rotational space-like
hypersurface of parabolic type.
\end{prop}
Combining with Proposition 3.12, we have
\begin{prop}
Any rotational space-like hypersurface of parabolic type in $\mathbb{R}^{n+1}_1$ is
a subprojective Riemannian manifold  with function $a+k^2 = 0$.
\end{prop}
The curvature tensor of a rotational space-like hypersurface of parabolic type has
the form
$$\mathcal R= -\frac{1}{R^2}\, \pi + \frac{R'^2}{R^2(RR''+R'^2)}\,\Phi.$$

\section{A local classification of Riemannian QC-manifolds}
Let $(M,g,\xi) \; (\dim M=n\geq 4)$ be a Riemannian QC-manifold.
Then the Riemannian curvature tensor $\mathcal R$ of $M$ has the form
$$ \mathcal R = a\pi +b \Phi.\leqno(4.1)$$

We consider manifolds free of points in which the tensor $\mathcal
R$ is of constant sectional curvatures, i.e. $b \neq 0$ in all
points of $M$.

We note that the condition $a=0$ implies that $b=0$.

In \cite {GM} we proved that a Riemannian QC-manifold with positive
horizontal sectional curvatures, i.e. $a > 0$, can be locally
embedded as a canal hypersurface in Euclidean space $\R^{n+1}$.

In this section we study Riemannian QC-manifolds with negative
horizontal sectional curvatures, i.e. $a < 0$.

The basic step in
our classification of Riemannian QC-manifolds is the following
theorem.
\begin{thm}
Let $(M,g,\xi) \; (\dim M=n\geq 4)$ be a Riemannian QC-manifold with curvature tensor
$(4.1)$ satisfying the conditions:
$$b\neq 0, \qquad a<0.$$
Then the manifold is locally $\xi$-isometric to a canal space-like hypersurface
in $\R^{n+1}_1$.

Moreover, the manifold is locally $\xi$-isometric to a canal space-like hypersurface
of elliptic, hyperbolic or parabolic type, according to
$$ a+k^2 >0, \quad a+k^2 <0 \quad {\rm or} \quad a+k^2 = 0,$$
respectively.
\end{thm}
{\it Proof:} Under the conditions of the theorem the curvature tensor of the manifold
$M$ has the form (4.1) and all equalities (2.3) - (2.7) are valid. We put
$$\alpha = \sqrt{-a}, \qquad \beta=-\frac{b}{\sqrt{-a}}$$
and consider the symmetric tensor
$$h=\sqrt{-a}\,g-\frac{b}{\sqrt{-a}}\,\eta \otimes \eta = \alpha \,g + \beta \,
\eta \otimes \eta \leqno(4.2)$$
on $M$.

An immediate verification shows that the curvature tensor $\mathcal R$ of
the manifold $(M,g,\xi)$ has the following construction
$$\mathcal R(X,Y,Z,U)=-\{h(Y,Z)\,h(X,U)-h(X,Z)\,h(Y,U)\},\leqno(4.3)$$
i.e.
$$\mathcal{R}=-(\alpha^2 \,\pi + \alpha \beta \,\Phi), \qquad a=-\alpha^2, \; b=-\alpha \beta.$$

We shall show that the tensor $h$ satisfies the Codazzi equation
$$(\nabla'_X h)(Y,Z)-(\nabla'_Y h)(X,Z)= 0, \quad X,Y \in \mathcal XM.\leqno(4.4)$$

Taking into account (4.2), we calculate
$$\begin{array}{ll}
(\nabla'_X h)(Y,Z)-(\nabla'_Y h) (X,Z)= & d \alpha(X) \, g(Y,Z)-d \alpha(Y)\, g(X,Z)\\
[2mm]
&+\,(d \beta(X)\, \eta(Y) - d \beta(Y)\, \eta(X))\,\eta(Z)\\
[2mm]
&+\,\beta\, d \eta(X,Y)\, \eta(Z)\\
[2mm]
&+\,\beta\,(\eta(Y)(\nabla'_X \eta)(Z)-\eta(X)(\nabla'_Y \eta)(Z)).
\end{array}\leqno(4.5)$$

We prove that the right hand side of (4.5) is identically zero. Since
any tangent vector  is decomposable as in (2.1), we divide the proof into four steps.
Taking into account that $a=-\alpha^2, \; b=-\alpha \beta$, we apply equalities
(2.3) - (2.7) and obtain consequently:

1) If $X=x,\, Y=y, \, Z=z $, then the right hand side of (4.5) reduces to
$$d \alpha(x) \, g(y,z) -d \alpha(y) \, g(x,z),$$
which is zero because of (2.3).

2) If $X=x, \, Y=y, \, Z= \xi$, then the right hand side of (4.5) reduces to
$$\beta \, d\eta(x,y),$$
which is zero in view of (2.4).

3) If $X=x, \, Y=\xi, \, Z= \xi$, then the right hand side of (4.5) reduces to
$$d\beta(x)+\beta\,d \eta(x,\xi),$$
which is zero as a consequence of (2.5).

4) If $X= \xi, \, Y=y, \, Z=z$, then the right hand side of (4.5) reduces to
$$ \xi(\alpha) \, g(y,z)-\beta (\nabla'_y \eta)(z),$$
which is zero because of (2.6).

Combining the above cases 1) - 4), we conclude that the right hand side of (4.5)
is equal to zero for all $X, Y, Z \in \mathcal XM $,
i.e. the tensor $h$ satisfies (4.4) identically.

Now we can apply the fundamental embedding theorem for hypersurfaces in $\R^{n+1}_1$
and obtain that the Riemannian QC-manifold $(M,g,\xi)$ can be locally embedded
as a hypersurface in the Minkowski space $\R^{n+1}_1$.

If $N$ is the unit normal vector field to a hypersurface with second fundamental
form $h$, then  the curvature tensor $\mathcal R$ of this hypersurface satisfies
the identity
$$\mathcal R(X,Y,Z,U)=g(N,N)\,\{h(Y,Z)\,h(X,U)-h(X,Z)\,h(Y,U)\}.$$
Comparing with (4.3) we obtain that the Riemannian QC-manifold $(M,g,\xi)$ is
embedded locally as a space-like hypersurface in $\R^{n+1}_1$. Further, we denote
this hypersurface again with $(M,g,\xi)$.

Now  $(M,g,\xi)$  is a space-like hypersurface in $\R^{n+1}_1$, whose second
fundamental form $h$ satisfies (4.2).

Next we prove that $M$ is locally a part of a space-like canal hypersurface
in $\R^{n+1}_1$.

Let $Z$ be the position vector field of $M$ and $p$ be a fixed point in $M$.
Denote by $S_p$ the maximal integral submanifold of the distribution $\Delta$
containing $p$. Using the property $d\alpha=\xi(\alpha) \,\eta$, we get
$\alpha = {\rm const}$ on $S_p$. Then the equality
$$\nabla'_xN=-\alpha \, x$$
implies that the vector function $Z-(1/\alpha)\,N$ is constant at the points
of $S_p$. We set
$$z=Z-\frac{1}{\alpha}\,N,$$
and conclude that $S_p$ lies on the time-like hypersphere $S^n$ with center $z$
and radius $R=(1/\alpha)$, and both hypersurfaces, $M$ and $S^n$, have the same
normals at the points of $S_p$.

Since the distribution $\Delta$ determines a one-parameter family of submanifolds
$Q^{n-1}(s), \, s \in J$ in a neighborhood $U$ of $p$, then $U$ is a part of the
envelope of this family.

Finally we apply Propositions 3.6, 3.9, 3.12 and obtain the second part of the theorem. \qed
\vskip 2mm
Applying Theorem 4.1, we obtain immediately
\begin{thm}
Let $(M,g,\xi) \; (\dim M=n\geq 4)$ be a subprojective Riemannian manifold with
curvature tensor $(4.1)$ satisfying the conditions:
$$b\neq 0, \qquad a<0.$$
Then the manifold is locally $\xi$-isometric to a rotational space-like hypersurface
in $\R^{n+1}_1$.

Moreover, the manifold is locally $\xi$-isometric to a rotational space-like
hypersurface of elliptic, hyperbolic or parabolic type, according to
$$ a+k^2 >0, \quad a+k^2 <0 \quad {\rm or} \quad a+k^2 = 0,$$
respectively.
\end{thm}

\section{Conformally flat Riemannian hypersurfaces in Euclidean or Minkowski space}

\subsection{Conformally flat hypersurfaces in Euclidean space}
A hypersurface $M$ in Euclidean space is said to be quasi-umbilical \cite{CY2} if its
second fundamental form $h$ satisfies the equality
$$h=\alpha \, g +\beta \, \eta \otimes \eta,\leqno(5.1)$$
for some functions $\alpha \neq 0$, $\beta \neq 0$, and a unit 1-form $\eta$ on $M$.
The close relation between conformally flat hypersurfaces in Euclidean
space from one hand side, and quasi-umbilical hypersurfaces in $\R^{n+1}$
from another hand side, is the following statement {\cite{C,Sch1} \,(see also \cite{NM}):
\begin{lem} {\bf (Cartan - Schouten)}
Let $M$ be a conformally flat hypersurface in $\R^{n+1}$. Then the shape operator
of $M$ at any point has a root of multiplicity at least $n-1$.
\end{lem}

As a result of Lemma 5.1 we have

\begin{lem}
Any conformally flat  hypersurface $M$ in Euclidean space, which is free of umbilical
points, is quasi-umbilical.
\end{lem}
{\it Proof:} Let $A$ be the shape operator of the hypersurface $M$. Since $M$ is
free of umbilical points, then according to Lemma 5.1 the operator $A$ has at any point
 two different eigenvalues $\alpha$ and $\alpha + \beta$ of multiplicity
$n-1$ and $1$, respectively. Let $\xi$ be the unit eigenvector field, corresponding
to the function $\alpha + \beta$. Denoting by $\eta$ the 1-form, corresponding to
$\xi$ with respect to the metric $g$, we obtain (5.1). \qed

Equality (5.1) implies that the curvature tensor $\mathcal R$ of $M$ has the form
$$\mathcal R= \alpha^2 \, \pi + \alpha \beta \, \Phi\,; \quad
\alpha^2 > 0, \quad \alpha \beta \neq 0, $$
i.e. $(M,g,\xi)$ is a Riemannian QC-manifold with positive horizontal sectional curvatures.

Applying Proposition 3 \cite{GM}, we obtain:

\begin{thm}
Any conformally flat hypersurface $M$ in $\R^{n+1}$, which is free of  umbilical
points, locally lies on a canal hypersurface.
\end{thm}

If $(M,g)$ is a conformally flat hypersurface in $\R^{n+1}$, then the manifold
$(M,g)$ admits a unit vector field $\xi$, such that $(M,g,\xi)$ is a Riemannian
QC-manifold with positive horizontal sectional curvatures $a > 0$. Any two
locally isometric conformally flat hypersurfaces are locally $\xi$-isometric,
i.e. rigid. Taking into account Theorem 2.4, we obtain that any isometric
embedding of a conformally flat Riemannian manifold into $\R^{n+1}$ is locally
determined up to a motion. We also recall the results of R. Beez \cite{B}
and W. Killing \cite{Kil}:

{\it A hypersurface in the Euclidean space is rigid if at least three principal
curvatures are different from zero at each point of it, i.e. the hypersurface
has type-number $\geq 3$ at each point.}

Thus, we obtained that
\vskip 1mm
\emph{The local theory of conformally flat Riemannian manifolds, isometrically
embedded as hypersurfaces in Euclidean space, is equivalent to the local theory
of Riemannian QC-manifolds with positive horizontal sectional curvatures.}
\vskip 1mm
Taking into account the local classification of hypersurfaces in Euclidean space of
constant sectional curvature and  Theorem 5.3, we obtain a
local geometric classification of conformally flat hypersurfaces in $\R^{n+1}$:

\begin{thm}
Any conformally flat hypersurface $M$ in $\R^{n+1}$ is locally a part of one of
the following hypersurfaces:

(i) \; hyperplane \; $(\alpha=\beta=0)$;

(ii) \; hypersphere \; $(\alpha \neq 0, \; \beta=0)$;

(iii) \; developable hypersurface \; $(\alpha = 0, \; \beta \neq 0)$;

(iv) \; canal hypersurface \; $(\alpha \neq 0, \; \beta \neq 0)$.
\end{thm}

\subsection{Conformally flat space-like hypersurfaces in Minkowski space}

Let $M$ be a space-like hypersurface in Minkowski space with second fundamental
form $h$. Similarly to the Euclidean case, we call the hypersurface $M$
{\it quasi-umbilical} if
$$h=\alpha \, g +\beta \, \eta \otimes \eta,\leqno(5.2)$$
for some functions $\alpha \neq 0$, $\beta \neq 0$, and a unit 1-form $\eta$ on $M$.

The proof of Lemma 5.1 is also valid without any essential changes for
conformally flat hypersurfaces in Minkowski space.

\begin{lem}
Let $M$ be a conformally flat space-like hypersurface in $\R^{n+1}_1$. Then the
shape operator of $M$ at any point has a root of multiplicity
 at least $n-1$.
\end{lem}

Analogously to Subsection 5.1., Lemma 5.5 implies the following statement.

\begin{lem}
Any conformally flat space-like hypersurface $M$ in Minkowski space, which is
free of  umbilical points, locally is quasi-umbilical.
\end{lem}

Equality (5.2) implies that the curvature tensor $\mathcal R$ of $M$ has the form
$$\mathcal R= -\alpha^2 \, \pi - \alpha \beta \, \Phi\,; \quad
\alpha^2 > 0, \quad \alpha \beta \neq 0, $$
i.e. $(M,g,\xi)$ is a Riemannian QC-manifold with negative horizontal sectional
curvatures.

Applying Theorem 4.1, we obtain:

\begin{thm}
Any conformally flat space-like hypersurface $M$ in Minkowski space, which is free
of  umbilical points, locally is a part of a canal space-like hypersurface.
\end{thm}

If $(M,g)$ is a conformally flat space-like hypersurface in $\R^{n+1}_1$, then the
manifold $(M,g)$ admits a unit vector field $\xi$, such that $(M,g,\xi)$ is a
Riemannian QC-manifold with negative horizontal sectional curvatures $a < 0$.
Any two locally isometric conformally flat space-like hypersurfaces are locally
$\xi$-isometric, i.e. rigid. All isometric realizations of a conformally
flat Riemannian manifold into $\R^{n+1}_1$ are locally congruent.

Thus we have:
\vskip 1mm
\emph{The local theory of conformally flat Riemannian manifolds,
isometrically immersed as space-like hypersurfaces in Minkowski space, is
equivalent to the local theory of Riemannian QC-manifolds with negative horizontal
sectional curvatures.}
\vskip 1mm
Taking into account the local classification of space-like hypersurfaces of
constant sectional curvatures in Minkowski space, and Theorem 5.7, we obtain
the following geometric classification of conformally flat space-like
hypersurfaces in Minkowski space:

\begin{thm}
Any conformally flat space-like hypersurface $M$ in Minkowski space is locally
a part of one of the following hypersurfaces:

(i) \; a space-like hyperplane \; $(\alpha=\beta=0)$;

(ii) \; a space-like hypersphere \; $(\alpha \neq 0, \; \beta=0)$;

(iii) \; a space-like developable hypersurface \; $(\alpha = 0, \; \beta \neq 0)$;

(iv) \; a space-like canal hypersurface \; $(\alpha \neq 0, \; \beta \neq 0)$.
\end{thm}

\end{document}